\documentclass[notitlepage,leqno,10pt]{article}

\usepackage[latin1]{inputenc}
\usepackage{amsfonts}
\usepackage{amsmath}
\usepackage{amssymb}
\usepackage{amsrefs}
\newcommand{\fint}{-\!\!\!\!\!\!\!\int}

\usepackage{color}
\definecolor{green2}{rgb}{0,0.6,0}

\newcommand{\R}{\mathbb{R}}

\newcommand{\proof}[1]{\par\smallskip\noindent{\bf Proof#1.}}
\newcommand{\qed}{\penalty 500\hfill$\square$\par\medskip}

\newcommand{\eps}{\varepsilon}
\newcommand{\pd}[2]{\frac{\partial #1}{\partial #2}}

\newcommand{\superp}{\overline{p_{0}} }
\newcommand{\subp}{\underline{p_{0}} }

\def\XXint#1#2#3{{\setbox0=\hbox{$#1{#2#3}{\int}$}
     \vcenter{\hbox{$#2#3$}}\kern-.5\wd0}}

\newtheorem{lem}         {Lemma}[section]
\newtheorem{pro}    [lem]{Proposition}
\newtheorem{thm}    [lem]{Theorem}
\newtheorem{rem}    [lem]{Remark}
\newtheorem{cor}    [lem]{Corollary}

\title{Stable solutions of $-\Delta u = f(u)$ in $\R^N$}
\date{\today}
\author{L. Dupaigne$^{1,2}$ and A. Farina$^{1,3}$}
\begin{document}
\maketitle

{\begin{center}
$^{1}${\small
LAMFA, UMR CNRS 6140, Universit\'e Picardie Jules Verne \\33, rue St Leu, 80039 Amiens, France \\
\smallskip
$^2$corresponding author, \tt louis.dupaigne@math.cnrs.fr\\
\smallskip
$^3$
alberto.farina@u-picardie.fr\\}
\end{center}}

\abstract{Several Liouville-type theorems are presented for stable solutions of the equation $-\Delta u = f(u)$ in $\R^N$, where $f>0$ is a general convex, nondecreasing functions. Extensions to solutions which are merely stable outside a compact set are discussed.}

\section{Introduction}
For $N\ge1$ and $f\in C^1(\R)$ consider the equation  
\begin{equation}\label{main}
-\Delta u = f(u)\qquad\text{in $\R^N$.}
\end{equation}
The aim of this paper is to classify solutions $u\in C^2(\R^N)$ which are stable i.e. such that  for all $\varphi\in C^1_{c}(\R^N)$,
\begin{equation}
\label{stability}
\int_{\R^N}f'(u)\varphi^2\; dx \le \int_{\R^N} |\nabla \varphi |^2\;dx.
\end{equation}
For some of our results, we shall assume in addition  $u>0$ in $\R^N$ and/or $u\in L^{\infty}(\R^N)$. We shall also discuss extensions to solutions which are merely stable outside a compact set (i.e. \eqref{stability} holds for test functions supported in the complement of a given compact set $K\subset\subset\R^N$).  

Stable \it radial \rm solutions of \eqref{main} are by now well-understood :  by the work of Cabr\'e and Capella \cite{cabre-capella}, refined by Villegas in \cite{villegas}, every bounded radial stable solution of \eqref{main} must be constant if $N\le10$. The result holds for any nonlinearity $f\in C^1(\R)$. Conversely, there exist unbounded radial stable solutions in any dimension. Take for example, $u(x)= \left| x \right|^2/2N$ solving \eqref{main} with $f(u)=-1$.  Also, there are examples of bounded radial stable solutions when $N\ge11$. See e.g. \cite{villegas}, \cite{farina-lane-emden}. When dealing with nonradial solutions, much less is known.  In the case $N=2$, any stable solution of \eqref{main} with bounded gradient is one-dimensional (i.e. up to a rotation of space, $u$ depends only on one variable) under the sole assumption that $f$ is locally Lipschitz continuous (see \cite{farina-valdinoci-sciunzi}). 
In arbitrary dimension, a complete analysis of stable solutions and solutions which are stable outside a compact set is provided for two important nonlinearities $f(u)=\left| u \right|^{p-1}u$, $p>1$ and $f(u)=e^{u}$ in \cite{farina-cras}, \cite{farina-lane-emden}, \cite{farina-cras2} and \cite{dancer-farina}. 

Under a mere nonnegativity assumption on the nonlinearity, we begin this paper by stating that up to space dimension $N=4$, bounded stable solutions of \eqref{main} are trivial : 
\begin{thm}\label{th1} 
Assume $f\in C^1(\R)$, $f\ge0$ and $1\le N\le4$. Assume $u\in C^2(\R^N)$ is a bounded, stable solution of \eqref{main}.
Then, $u$ is constant.
\end{thm}
\begin{rem}
It would be interesting to know whether Theorem \ref{th1} still holds if one assumes that $u$ is unbounded but $\nabla u$ is bounded. 
\end{rem}
\subsection{Power-type nonlinearities}
For our next set of results, we restrict to the following class of nonlinearities
\begin{equation} \label{convex class} 
f\in C^0(\R^+)\cap C^2(\R^+_{*}), \text{$f>0$ is nondecreasing and convex in $\R^+_{*}$}.
\end{equation} 
As demonstrated in \cite{farina-lane-emden} for the particular case of the power nonlinearities $f(u)=\left| u  \right|^{p-1}u$, two critical exponents play an important role, namely the classical Sobolev exponent
\begin{equation} \label{sobolev exponent}
p_{S}(N) = \frac{N+2}{N-2},\qquad\text{for $N\ge3$} 
\end{equation} 
and the Joseph-Lundgren exponent
\begin{equation} \label{pcofn} 
 p_{c}(N)=\frac{(N-2)^2-4N+8\sqrt{N-1}}{(N-2)(N-10)},\qquad\text{for $N\ge11$.} 
\end{equation} 
In order to relate the nonlinearity $f$ and the above exponents, we introduce a quantity $q$ defined  for $u\in\R^+_{*}$ by 
\begin{equation}\label{q of u}
q(u)=\frac{f'^2}{ff''}(u)=\frac{(\ln f)'}{(\ln f')'}(u)
\end{equation}
whenever  $ff''(u)\neq 0$, $q(u)=+\infty$ otherwise.
When $f(u)=\left| u \right|^{p-1}u$, $p\ge1$, $q$ is independent of $u$ and coincides with the conjugate exponent of $p$ i.e. $\frac{1}{p}+ \frac{1}{q} =1$.
In this section, we assume that $q(u)$ converges as $u\to0^+$ and denote its limit :
\begin{equation}
\label{equation q}
q_{0} = \lim_{u\to 0^+} q(u)\in \overline{\R}.
 \end{equation}
\begin{rem}If $u\in C^2(\R^N)$, $u\ge0$ solves \eqref{main}, and  \eqref{convex class} holds, then $f(0)=0$. 
\end{rem}
In dimension $N=1,2$, this follows directly from the classical Liouville theorem for superharmonic nonnegative functions. For a proof in dimension $N\ge3$, see Step 6. in  Section 6. We then observe that 
\begin{lem}\label{basic lemma} 
If $f\in C^0(\R^+)\cap C^2(\R^+_{*})$ is convex nondecreasing, $f(0)=0$ and \eqref{equation q} holds, then in fact $q_{0}\in[1,+\infty]$.
\end{lem} 
\proof{} 
Indeed, assume by contradiction there exists $\theta>1$ such that $0\le q(u)\le 1/\theta$ in a neighbourhood of $0$. Consequently,  near $0$,
$$
\frac{f''}{f'} - \theta\frac{f'}{f}\ge 0. 
$$
So, ${f'}/{f^\theta}$ is nondecreasing hence bounded above near $0$. Integrating again, we deduce that $f^{1-\theta}(u)\le Cu + C'$ near $0$, which is not possible if $f(0)=0$. \qed 
\noindent Define now $p_{0}\in \overline \R$, the conjugate exponent of $q_{0}$ by
\begin{equation} \label{pnought}
1/p_{0}+1/q_{0}=1.
\end{equation}  
The exponent $p_{0}$ must be understood as a measure of the ``flatness`` of $f$ at $0$. 
All nonlinearities $f$ such that \eqref{convex class} holds and which either are analytic at the origin or have at least one non-zero derivative at the origin or are merely  of the form $f(u)=u^p g(u)$, where $p\ge 1$ and $g(0)\neq0$, satisfy \eqref{equation q}. Exponentially flat functions such as $f(u)=e^{-1/u^2}$ also qualify (with $p_{0}=+\infty$).  However, there should exist (convex increasing) nonlinearities failing \eqref{equation q}. This being said, we establish the following theorem.
\begin{thm}\label{th2}
Assume $f\in  C^0(\R^+)\cap C^2(\R^+_{*})$ is nondecreasing, convex, $f>0$ in $\R_*^+$ 
and \eqref{equation q} holds.
Assume $u\in C^2(\R^N)$ is a bounded,nonnegative, stable solution of \eqref{main}.
Then, $u\equiv0$ if either of the following conditions holds
\begin{enumerate}
\item $1\le N\le 9$,
\item $N=10$ and $p_{0}<+\infty$, where $p_{0}$ is given by \eqref{pnought},
\item $N\ge11$ and $p_{0}<p_{c}(N)$, where $p_{0}$ is given by \eqref{pnought}  and $p_{c}(N)$ by \eqref{pcofn}  
\end{enumerate}
\end{thm}

\begin{rem}\label{optimality} 
Theorem \ref{th2} was first proved by A. Farina,when $f(u)=\left| u  \right|^{p-1}u$. See \cite{farina-lane-emden}.
As observed e.g. in \cite{farina-lane-emden}, for $N\ge11$, there exists a non constant bounded positive stable solution for $f(u)=\left| u  \right|^{p-1}u$ as soon as $p\ge p_{c}(N)$. So our result is sharp in the class of power-type nonlinearities for $N\ge11$. We do not know whether Theorem \ref{th2} remains true when $N=10$ and $p_{0}=+\infty$. 
We do not know either if for $N\le10$, assumption \eqref{equation q} can be completely removed. See Theorem \ref{thbeyond} in Section \ref{beyond} for partial results in this direction.  
See also \cite{villegas} for a positive answer in the radial case. 
\end{rem}


\subsection{Some generalizations : unbounded and sign-changing solutions, beyond power-type nonlinearities}\label{beyond} 
First, we discuss the case of unbounded solutions.
When $f(u)=\left| u  \right|^{p-1}u$, the assumption $u\in L^\infty(\R^N)$ is unnecessary, see \cite{farina-lane-emden}. For general power-type nonlinearities, Theorem \ref{th2} remains true for unbounded solutions under an  additional assumption on the behaviour of $f$ at $+\infty$ : 
\begin{cor}\label{unbounded} 
Assume as before that $f\in  C^0(\R^+)\cap C^2(\R^+_{*})$ is nondecreasing, convex, $f>0$ in $\R_*^+$ 
and \eqref{equation q} holds.
Let $\underline{p_{\infty}}\in\overline\R$ defined by
\begin{align}
\label{equation q infty}
\overline{q_{\infty}} := \limsup_{u\to +\infty} q(u)&,\\
1/\underline{p_{\infty}} + 1/\overline{q_{\infty}} =1.& \nonumber
\end{align}
Let $u\in C^2(\R^N)$ denote a nonnegative, stable solution of \eqref{main}.
Then, $u\equiv0$ if either of the following conditions hold
\begin{enumerate}
\item $1\le N\le 9$ and 
$1<\underline{p_{\infty}}$,
\item $N=10$, $p_{0}<+\infty$ and  $1<\underline{p_{\infty}}<+\infty$,
\item $N\ge11$, $p_{0}<p_{c}(N)$ and $1<\underline{p_{\infty}}<p_{c}(N)$.\end{enumerate}
\end{cor}
Next, we look at solutions which may change sign.
When $f(u)=\left| u  \right|^{p-1}u$, the assumption $u\ge0$ is also unnecessary, see \cite{farina-lane-emden}. For power-type nonlinearities, Theorem \ref{th2} can be extended to the case of solutions of arbitrary sign if $f$ is odd :
\begin{cor}\label{sign changing}
Assume that $f\in  C^0(\R)\cap C^2(\R^+_{*})$ is nondecreasing and that when restricted to $\R_*^+$, $f$ is convex and $f>0$.
Assume \eqref{equation q} holds. Assume in addition that $f$ is odd. Let $u\in C^2(\R^N)$ denote a bounded, stable solution of \eqref{main}.
Then, $u\equiv0$ if either of the following conditions hold
\begin{enumerate}
\item $1\le N\le 9$ and $1<p_{0}$,
\item $N=10$ and $1<p_{0}<+\infty$,
\item $N\ge11$ and $1<p_{0}<p_{c}(N)$.
\end{enumerate}
\end{cor} 

\begin{rem}\label{remark odd}The above Corollary remains true if $f$ is not odd but simply if $f(0)=0$ and the assumptions made on $f$ also hold for $\tilde f$ defined for $u\in\R^+$ by $\tilde f(u)=- f (-u)$.
\end{rem}

\begin{cor}\label{sign changing 2} Assuming in addition  $1<\underline{p_{\infty}}$ if $N\le 9$ (respectively $1<\underline{p_{\infty}}<+\infty$ if $N=10$ and $1<\underline{p_{\infty}}<p_c(N)$ when $N\ge11$), Corollary \ref{sign changing}   remains valid for any stable solution. That is, one can drop the assumptions $u\ge0$ and $u\in L^\infty(\R^N)$.  
\end{cor}
Finally, we study nonlinearities for which \eqref{equation q} fails. To do so, we introduce 
$\overline{q_{0}}, \underline{q_{0}}\in\overline\R$ defined by
 \begin{equation}
\label{equation q sup}
\overline{q_{0}} = \limsup_{u\to 0^+} q(u), \quad \underline{q_{0}} = \liminf_{u\to 0^+} q(u).
 \end{equation}

\begin{thm}\label{thbeyond}
Assume $f\in C^0(\R^+)\cap C^2(\R^+_{*})$ is nondecreasing, convex, $f>0$ in $\R_*^{+}$  
and let 
$\overline{q_{0}}, \underline{q_{0}}$ defined by \eqref{equation q sup}. 
Assume $u\in C^2(\R^N)$ is a bounded, nonnegative, stable solution of \eqref{main}.
Then, $u\equiv0$ if either of the following conditions hold
\begin{enumerate}
\item $3\le N$ and $\underline{q_{0}}>\frac N{2}$,\label{case 1} 
\item $1\le N\le 6$ and $\overline{q_{0}}<\infty$,\label{case 2} 
\item \label{case 3} 
$1\le N$ and $\frac{4}{N-2}\left(1+1/\sqrt{\overline{q_{0}}}\right)>1/\underline{q_{0}}$.
\end{enumerate}
\end{thm}

\begin{rem}
The above theorem is of particular interest when $f'$ is convex or concave near the origin. Assume $f(0)=f'(0)=0$ (this is not restrictive, see Remark \ref{remark pohozaev}). Apply Cauchy's  mean value theorem : given $u_{n}\in\R^+_{*}$, there exists $v_n\in(0,u_{n})$ such that
$$
q(u_{n}) = \left.\frac{f'^2}{ff''}\right|_{u=u_{n}} = \left. \frac{2f'f''}{f'f''+ff'''}\right|_{u=v_{n}}.   
$$
If $f'''\ge0$ near $0$, we deduce that $\overline{q_{0}}\le2$. By case \ref{case 2} of the Theorem, we conclude that if $f'$ is convex near $0$ and $N\le6$, then $u\equiv0$. Similarly, if $f'$ is concave near $0$, $\underline{q_{0}}\ge2$. By case \ref{case 3} of the Theorem, we conclude that  if $f'$ is concave near $0$ and $N\le9$ (or $N=10$ and $\overline{q_{0}}<+\infty$), then $u\equiv0$.    
\end{rem}
\begin{rem}
Our methods yield absolutely no result under the assumption $10\ge N\ge5$ and
$$
\underline{q_{0}}\le \frac {N-2}{4}<\overline{q_{0}}=\infty.
$$
\end{rem}

\subsection{Solutions which are stable outside a compact set}
Set aside the case where $f$ is a power or an exponential nonlinearity, little is known about the classification of solutions of \eqref{main} which are stable outside a compact set. Even in the radial case. Now, recall the definition of the critical exponents given in \eqref{sobolev exponent} and \eqref{pcofn}. As demonstrated in \cite{farina-lane-emden}, the nonlinearities $f(u)=\left| u  \right|^{p-1}u$, $p=p_{S}(N)$, $N\ge3$ and $p\ge p_c(N)$, $N\ge11$ must be singled out. For such values of $p$,  radial solutions which are stable outside a compact set are nontrivial and completely classified, while for other values of $p>1$, all solutions which are stable outside a compact set (whether radial or not) must be constant. See \cite{farina-lane-emden}.
When dealing with more general nonlinearities, the first basic step consists in determining the behaviour of a solution $u$ at infinity. This can be done by exploiting the classification of stable solutions obtained in Theorem \ref{th2} and Corollary \ref{sign changing}  :  
\begin{pro}\label{prop1}  Assume  $f\in C^0(\R)$.
Assume $u=0$ is the only bounded stable $C^2$ solution of \eqref{main}. If $u\in C^2(\R^N)$ is a bounded solution 
of \eqref{main} which is stable outside a compact set, then,
$$
\lim_{\left| x\right| \to \infty } u(x) =0.
$$  
\end{pro} 
\begin{rem}
As follows from the proof, the same result is valid for bounded \rm positive \it solutions which are stable outside a compact set, under the weaker assumption that all bounded \rm positive \it stable solutions of the equation are constant.
\end{rem}
\begin{rem} \label{no solution} 
If $f'(0)>0$, then in fact there exists no bounded solution of \eqref{main} which is stable outside a compact set. See the proof of Proposition \ref{prop1}.  
\end{rem} 
\begin{rem} \label{zeros of f} 
Clearly, if we assume instead that $f$ vanishes only at $u_{0}\neq 0$ then  $\lim_{\left| x\right| \to \infty } u(x) =u_{0}$. Similarly, we leave the reader check that if the set of zeros of $f$ is totally disconnected and the only bounded stable solutions of the equation are constant, then $\lim_{\left| x\right| \to \infty } u(x) =u_{0}$, where $u_{0}$ is a zero of $f$.
\end{rem}  
\begin{rem}\label{unbounded prop}
We do not know if a version of Proposition \ref{prop1} holds if one assumes that $f$ vanishes only at $-\infty$ or $+\infty$. If $f(u)=e^{u}$ and $N=2$ (see e.g. \cite{farina-cras2}), there exist (infinitely many) solutions of \eqref{main} which are stable outside a compact set and such that 
$$
\lim_{\left| x\right| \to \infty } u(x) =-\infty.
$$ 
\end{rem}   

\proof{ of Proposition \ref{prop1}} For $k\ge1$, let $\tau_{k}\in\R^N$  such that $\lim_{k\to\infty}| \tau_{k}| = +\infty$ and let $u_{k}(x)=u(x+\tau_{k})$ for $x\in\R^N$. Standard elliptic regularity implies that a subsequence of $(u_{k})$ converges in the topology of $C^2_{\rm loc}(\R^N)$ to a solution $v$ of \eqref{main}. In addition, since $u$ is stable outside a compact set, $v$ is stable. Therefore, $v$ is constant and $f(v)=0$, so $v=0$. If $f'(0)>0$, then $v=0$ is clearly unstable, which is absurd. This proves Remark \ref{no solution}.
In addition, since $v=0$ is the unique cluster point of $(u_{k})$, the whole sequence must converge to $0$. Proposition \ref{prop1} follows.\qed

In light of Proposition \ref{prop1}, it is natural to try to characterize the speed of decay of our solutions as $\left| x\right| \to\infty$.  When $f$ is power-type, we have the following:

\begin{thm}\label{th3}
Assume $f\in C^0(\R^+)\cap C^2(\R^+_{*})$ is nondecreasing, convex, $f>0$ in $\R^+_*$, $f(0)=0$ 
and \eqref{equation q} holds.
Assume $u\in C^2(\R^N)$ is a bounded positive solution of \eqref{main}, which is stable outside a compact set. If either of the following conditions holds
\begin{enumerate}
\item $1\le N\le 9$,
\item $N=10$ and $p_{0}<+\infty$,
\item $N\ge11$ and $p_{0}<p_{c}(N)$, 
\end{enumerate}
then there exists a constant $C>0$ such that for all $x\in\R^N$ sufficiently large,
\begin{equation} \label{asymptotics}
u(x)\le C s(|x|).
\end{equation}   
In the above inequality, the speed of decay $s(R)$ is defined for $R>0$ as the unique solution $s=s(R)$ of  
\begin{equation} \label{speed def} 
f(A_{1}R^2f(s)) =A_{2}f(s),
\end{equation} 
where $A_{1},A_{2}$ are two positive constants depending on $N$ only. In other words, $s$ is given by $s(R)=f^{-1}\left(C_{1}R^{-2}g(C_{2}R^{-2})\right)$ where $C_{1},C_{2}$ are two positive constants depending on $N$ only and $g$ is the inverse function of $t\mapsto f(t)/t$.
\end{thm}
\begin{rem}\label{f sur t est inversible} 
In the above theorem, we have implicitly assumed that the functions $f$ and $t\to f(t)/t$ are invertible in a neighborhood of $0$. This is indeed true : by convexity of $f$, $t\to f(t)/t$ is nondecreasing. By Step 6 in Section 6, we must have $f(0)=0$ and $\lim_{t\to0^+}\frac{f(t)}{t}=0$. If there existed two values $0<t_{1}<t_{2}$ such that $\frac{f(t_{1})}{t_{1}}=\frac{f(t_{2})}{t_{2}}$, then, by convexity, $f$ would be linear on $(t_{1},t_{2})$, hence on $(0,t_{2})$ by convexity. This contradicts  $\lim_{t\to0^+}\frac{f(t)}{t}=0$. So, $t\to f(t)/t$ is invertible for $t>0$ small and so must be $f$. 
\end{rem}

\begin{rem}
Equation \eqref{speed def} looks somewhat complicated at first glance. For many nonlinearities (including $f(u)= \left| u \right|^{p-1}u$), one can actually set the constants $A_{1}, A_{2}$ equal to $1$. \eqref{speed def} then takes the simplified form
$$
\frac{f(s)}{s} = R^{-2}.
$$ 
In particular, when $f(u)= \left| u \right|^{p-1}u$, we recover the familiar speed $s(R)=R^{-\frac2{p-1}}$.
\end{rem}
\begin{rem} \label{estimate of h} If $p_{0}<\infty$, for all $\epsilon>0$, there exists $C>0$ such that 
$$
s(R)\le C R^{-\frac2{(p_{0}-1)}+\epsilon}\qquad\text{for $R\ge1$}.
$$
However, even when $p_{0}<\infty$, there should exist nonlinearities $f$ failing the estimate $s(R)\le C R^{-2/(p_{0}-1)}$.
\end{rem}  

\proof{ of Remark \ref{estimate of h}} 
An easy calculation shows that for all $\delta>0$ small, there exists  $C,\eps>0$ such that $C^{-1}u^{p_{0}+\delta}\le f(u)\le C u^{p_{0}-\delta}$ and $C^{-1}u^{p_{0}+\delta-1}\le f'(u)\le C u^{p_{0}-\delta-1}$ for $u\in(0,\eps)$ provided \eqref{equation q} holds and $p_{0}<+\infty$. Plugging this information into the definition of $s(R)$ yields  the desired conclusion.
\qed


From here on, our aim is to prove a Liouville-type result for solutions which are stable outside a compact set. As follows from the analysis in \cite{farina-lane-emden}, we must distinguish the sub and the supercritical case. We first consider the case where $p_{0}$ is subcritical i.e. 
\begin{equation}\label{subcritical range}  
p_{0}<\infty, \;\; N\le2\quad\text{ or }\quad p_{0}< p_{S}(N), \;\; N\ge3.
\end{equation} 
In this case, we make the following extra global assumption on $f$ :
\begin{equation} 
 \label{global assumption subcritical}  
(p_{0}+1)F(s)\ge sf(s)\qquad\text{for all $s\in\R$},
\end{equation}
where $F$ denotes the antiderivative of $f$ vanishing at $0$. Then, we have

\begin{thm} \label{thsubcritical}
Assume $f\in C^0(\R^+)\cap C^2(\R^+_{*})$ is nondecreasing, convex, $f>0$ in $\R^*_+$ 
and \eqref{equation q} holds.
Assume $u\in C^2(\R^N)$ is a bounded, nonnegative solution of \eqref{main}, which is stable outside a compact set. Assume
$p_{0}$ is subcritical (i.e. \eqref{subcritical range} holds) and $f$ satisfies the global inequality \eqref{global assumption subcritical}. Then, $u=0$. 
\end{thm}
We turn next to the supercritical case. We say that $p_{0}$ is in the supercritical range if 
 \begin{equation}\label{supercritical range} 
 \left\{
 \begin{aligned} 
 p_{S}(N)<p_{0}&<+\infty,&\;\; 3\le N\le10,\\ 
 &\text{ or }&\\
 p_{S}(N)<p_{0}&< p_{c}(N),&\;\; N\ge11.
\end{aligned}
\right. 
 \end{equation} 
 In this case, we begin by showing that the  asymptotic decay estimate \eqref{asymptotics} can be further improved. Namely, we show that not only $u(x)=O(s(\left| x \right|  ))$ but in fact  $u(x)= o(s(\left| x \right|  ))$. The price we pay is the following set of assumptions : we request that \it near the origin\rm, there exist constants $\eps, c_{1}, c_{2}>0$ such that
\begin{align}
\label{lower bound} 
f(u)&\ge c_{1}u^{p_{0}}\qquad\text{for $u\in(0,\eps)$}\\
\label{upper bound}
f'(u)&\le c_{2}u^{p_{0}-1}\qquad\text{for $u\in(0,\eps)$}.
\end{align}
By convexity of $f$, the above inequalities reduce to one when $f(0)=0$ : 
\begin{equation} \label{lower and upper bound} 
c_{2}u^{p_{0}}\ge uf'(u)\ge f(u)\ge c_{1}u^{p_{0}},\qquad\text{for $u\in(0,\eps)$}.
\end{equation} 
Compare this assumption with the already known estimate given in the proof of Remark \ref{estimate of h}.    
\begin{thm}\label{cor fine asymptotics}
Make the same assumptions as in Theorem \ref{th3}. Assume in addition that $f$satisfies the local estimates \eqref{lower bound}, \eqref{upper bound}. For $p_{0}$ in the supercritical range \eqref{supercritical range}, any bounded positive solution $u\in C^2(\R^N)$ of \eqref{main}, which is stable outside a compact set, satisfies
\begin{equation} \label{fine asymptotics}
u(x) = o\left(\left|  x \right|^{-\frac2{p_{0}-1}}\right) \quad\text{ and }\left| \nabla u(x) \right|  = o\left(\left|  x \right|^{-\frac2{p_{0}-1} -1}\right)\quad\text{as $\left| x \right| \to\infty$.}
\end{equation}     
\end{thm}
Finally, to obtain the Liouville theorem in the supercritical range, we assume in addition that
 \begin{equation} 
\label{global assumption}
 (p_{0}+1)F(s)\le sf(s)\qquad\text{for all $s\in\R$.}
 \end{equation} 
Note that the inequality is reversed compared to \eqref{global assumption subcritical}. 
Also note that since $f$ is nondecreasing, we automatically have $F(s)\le sf(s)$. \eqref{global assumption} can thus be seen as an improved global convexity assumption on $F$. 
We have
\begin{thm} \label{th4}
Assume $f\in C^2(\R^+)$ is nondecreasing, convex, $f>0$ in $\R^*_+$ 
and \eqref{equation q} holds.
Assume $u\in C^2(\R^N)$ is a bounded, nonnegative solution of \eqref{main}, which is stable outside a compact set. Assume
$p_{0}$ is in the supercritical range \eqref{supercritical range} and $f$ satisfies the local bounds \eqref{lower bound}, \eqref{upper bound} as well as the global inequality \eqref{global assumption}.  Then, $u\equiv0$.   
\end{thm}
\begin{rem}
As mentioned in Remark \ref{optimality}  , the above theorem is false for exponents $p_{0}\ge p_{c}(N)$, $N\ge 11$ or $p_{0}=p_{S}(N)$, $N\ge3$. \end{rem}

\begin{rem}For the nonlinearity $f(u) = \left| u \right|^{p-1}u$, \it all \rm the extra assumptions \eqref{global assumption subcritical}, \eqref{global assumption},  \eqref{lower bound}, \eqref{upper bound} are automatically satisfied.    
\end{rem}

\noindent The rest of the paper is organized as follows. In Section 2, we prove Theorem \ref{th1}.  Theorem \ref{th2} is the object of Section 3. In Section 4, we discuss the extensions given in Corollaries \ref{unbounded}, \ref{sign changing} and \ref{sign changing 2}. Theorem \ref{thbeyond}, which deals with nonlinearities which are not of power-type, is proved in Section 5. Section 6 is devoted to the proof of Theorem \ref{th3}, pertaining to the rate of decay of solutions which are stable outisde a compact set. The refined asymptotics obtained in Corollary \ref{cor fine asymptotics} is also derived in this section. Section 7 covers Theorem \ref{thsubcritical}, dealing with subcritical nonlinearities, while the supercritical case is addressed in Section 8.        
\section{The case of low dimensions $1\le N\le 4$ : proof of Theorem 1.1}
The proof bears resemblences with an argument found in \cite{ambrosio-cabre}. It relies on two simple arguments : a growth estimate of the Dirichlet energy on balls and a Liouville-type result for certain divergence-form equations (mainly due to Berestycki, Caffarelli and Nirenberg \cite{berestycki-caffarelli-nirenberg}), which applies to solutions with controlled energy. The specific form of the afore-mentioned equation is obtained by linearizing \eqref{main} and taking advantage of the stability assumption. The limitation $N\le4$ arises from the energy estimate on balls.
\proof{} For $R>0$, let $B_{R}$ denote the ball of radius $R$ centered at the origin. We begin by proving that there exists a constant $C>0$ independent of $R>0$ such that
\begin{equation}\label{estimate}
    \int_{B_R} |\nabla u|^2 \;dx \le C R^{N-2}.
\end{equation}
Let $M\ge\|u\|_\infty$, $\varphi\in C^2_c(\R^N)$ and multiply \eqref{main} by $(u-M)\varphi$ :
$$
\int_{\R^N} -\Delta u (u-M)\varphi \;dx = \int_{\R^N} f(u) (u-M)\varphi \;dx.
$$
Integrating by parts and recalling that $f\ge0$, it follows that
\begin{align*}
\int_{\R^N} |\nabla u|^2\varphi\;dx + \int_{\R^N} (u-M)\nabla u \nabla \varphi\;dx &= \int_{\R^N} f(u) (u-M)\varphi \;dx\\&\le 0,
\end{align*}
whence,
\begin{align*}
\int_{\R^N} |\nabla u|^2 \varphi\;dx &\le -\int_{\R^N} \frac{1}{2}\nabla(u-M)^2\nabla\varphi\;dx = \int_{\R^N} \frac{(u-M)^2}{2} \Delta\varphi\;dx\\
&\le 2M^2 \int_{\R^N}\left| \Delta\varphi \right|\;dx .
\end{align*}
Let  $\varphi_0$ denote any nonnegative test function such that $\varphi_0=1$ on $B_1$ and apply the above inequality with $\varphi(x)=\varphi_0(x/R)$. We obtain \eqref{estimate}.

Since $u$ is stable, there exists a solution $v>0$ of the linearized equation
\begin{equation}\label{linearized equation}
    -\Delta v = f'(u)v\qquad\text{in $\R^N$}.
\end{equation}
Let $\sigma_j = \frac1v \frac{\partial u}{\partial x_j}$ for $j=1,\dots,N$. Then, since $v$ and ${\partial u}/{\partial x_j}$ both solve the linearized equation \eqref{linearized equation}, it follows that 
\begin{equation}\label{divergence form equation}
    -\nabla\cdot\left(v^2\nabla \sigma_j\right)=0\qquad\text{in $\R^N$}.
\end{equation}
It is known that any solution $\sigma\in H^1_{\rm loc \it}(\R^N)$ of \eqref{divergence form equation} such that  
\begin{align*}
    \int_{B_R} v^2\sigma^2 \le CR^2,
\end{align*}
must be constant (see Proposition 2.1 in \cite{ambrosio-cabre}).
By \eqref{estimate}, we deduce that if $N\le4$, then  $\sigma_j$ is constant, i.e. there exists a constant $C_j$ such that
$$
\pd{u}{x_j}=C_j v.
$$
In particular, the gradient of $u$ points in a fixed direction i.e. $u$ is one-dimensional and solves 
\begin{align*}
    -u''=f(u)\qquad\text{in $\R$}.
\end{align*}
Since $f\ge0$ and $u$ is bounded, this is possible only if $u$ is constant and $f(u)=0$.

\hfill\qed

\section{The Liouville theorem for stable solutions : proof of Theorem \ref{th2}}
The proof is split into two separate cases, according to the value of $q_{0}$. We first consider the case $q_{0}>\frac{N}{2}$. It suffices to prove the following lemma. 
\begin{lem}\label{lemma2} Assume $f\in C^2(\R^+)$, $f>0$ is nondecreasing, convex and 
$$\underline{q_{0}}:=\liminf_{u\to 0^+}q(u)>N/2.$$
Assume $u\in C^2(\R^N)$ , $u\ge0$ and
\begin{equation} \label{diff eq}  
-\Delta u \ge f(u)\qquad\text{in $\R^N$}.
\end{equation} 
Then, $u\equiv0$.
\end{lem}
\begin{rem}
A stronger version of the above lemma has been recently proved by L. D'Ambrosio and E. Mitidieri (\cite{dambrosio-mitidieri}). 
\end{rem}

\proof{} 
Assume by contradiction that $u\neq0$. By the Strong Maximum Principle, $u>0$.
 
\noindent{\bf Step 1.}
Since $\underline{q_{0}}>\frac{N}{2}$, there exists $q>\frac{N}{2}$ such that 
\begin{equation*}
 \frac{f''f}{f'^2}<\frac{1}{q},
\end{equation*}
in a neighborhood of $0$. Equivalently, $\frac{f''}{f'}-\frac{1}{q}\frac{f'}{f}<0$. Hence, the function $\frac{f'}{f^{1/q}}$ is decreasing near $0$. In particular, there exists a constant $C>0$ such that
$
 \frac{f'}{f^{1/q}}\ge C
$
near $0$, which implies that for some $p<\frac{N}{N-2}$, $c_{1}>0$,
\begin{equation}
 \label{serrin}
f(u)\ge c_{1}u^{p}.
\end{equation}
The above inequality holds in a neighborhood of $0$. 

\noindent{\bf Step 2.} Since $p<p_{S}(N)$, there exists $\varphi>0$ solving
\begin{equation}\label{equation varphi} 
 \left\{
 \begin{aligned} 
 -\Delta \varphi &=c_{1}\varphi^p&\quad\text{in $B_{1}$}\\
 \varphi &=0&\quad\text{on $\partial B_{1}$}. 
\end{aligned}
\right. 
 \end{equation}  
We are going to prove that a rescaled version of $\varphi$ must lie below $u$. Let indeed $R>0$ and $\varphi_{R}(x)=R^{-\frac2{p-1}}\varphi(x/R)$ for $x\in B_{R}$, $\varphi_{R}(x)=0$ for $\left| x \right|\ge R$. Then,
\begin{equation*}
 \left\{
 \begin{aligned} 
 -\Delta \varphi_{R} &=c_{1}(\varphi_{R})^p&\quad\text{in $B_{R}$}\\
 \varphi_{R} &=0&\quad\text{on $\partial B_{R}$}. 
\end{aligned}
\right. 
 \end{equation*}  
Furthermore, since $p<\frac{N}{N-2}$,
\begin{equation} \label{scaling law} 
\frac{\| \varphi_{R} \|_{L^\infty(B_{R})} }{R^{2-N}}\le \frac{R^{-\frac2{p-1}}}{R^{2-N}}\| \varphi \|_{L^\infty(B_{1})}\to 0,\quad\text{as $R\to+\infty$}. 
\end{equation} 
\noindent{\bf Step 3.} Since $u>0$ is superharmonic, there exists a constant $c>0$ such that 
\begin{equation} \label{rate of decay from below} 
u(x)\ge c \left| x \right|^{2-N}\quad\text{for $\left| x \right| \ge 1$.} 
\end{equation} 
Indeed, the above inequality clearly holds for $\left| x \right| =1$, with $c=\min_{[\left| x \right| =1]}u$. In addition, the function $z=u - c \left| x \right|^{2-N}$ is superharmonic in
${[1\le\left| x \right| \le M]}$. By the Maximum Principle, $z\ge \min(0,\min_{{[\left| x \right| =M]}}z(x))$, in ${[1\le\left| x \right| \le M]}$. Hence, 
$z\ge \liminf_{M\to\infty}\min(0,\min_{{[\left| x \right| =M]}}z(x))=0$. \eqref{rate of decay from below} is established.

\noindent{\bf Step 4.} Collecting \eqref{scaling law} and \eqref{rate of decay from below}, we obtain for $R>0$ sufficiently large
$$
u\ge \varphi_{R}.
$$  
We conclude using the celebrated sliding method : first, by \eqref{scaling law}, $\| \varphi_{R} \|_{\infty}\to0$ as $R\to\infty$, so that by \eqref{serrin}, $f(\varphi_{R})\ge c_{1}(\varphi_{R})^p$, provided $R$ is sufficiently large. In particular,
$$
-\Delta (u-\varphi_{R})\ge f(u)-f(\varphi_{R})\ge0.
$$
By the Strong Maximum Principle, $u>\varphi_{R}$.
Next, we slide $\varphi_{R}$ in a given direction, say $\tilde \varphi_{R,t}(x)=\varphi_{R}(x+te_{1})$, where $e_{1}=(1,0,\dots,0)$. We want to prove that $u\ge \tilde \varphi_{R,t}$ for all $t\ge0$. If not, there exists $t_{0}\in(0,+\infty)$ such that $u\ge\tilde \varphi_{R,t_{0}}$ and $u(x_{0})=\tilde\varphi_{R,t_{0}}(x_{0})$ at some point $x_{0}\in\R^N$. But again we have
$$
-\Delta (u-\varphi_{R,t_{0}})\ge f(u)-f(\tilde\varphi_{R,t_{0}})\ge0.
$$   
and the Strong Maximum Principle would imply that $u\equiv \tilde \varphi_{R,t_{0}}$. This is not possible since $\varphi_{R,t_{0}}$ is compactly supported while $u$ is not. The above argument holds if $e_{1}$ is replaced by any other direction $e\in S^{N-1}$. In particular, $u\ge\max{\varphi_{R}}>0$, which is possible, since $u$ is superharmonic, only if $u$ is constant. Since, $f>0$, we obtain a contradiction. Hence, $u\equiv0$.

 \;\hfill\qed
\begin{rem}\label{remark pohozaev}
If $f(0)\neq0$ or $f'(0)\neq0$, then \eqref{serrin} clearly holds in a neighborhood of  $0$ and we may work as above to conclude that $u$ is constant. We may therefore assume for the rest of the proof that $f(0)=f'(0)=0$.
\end{rem}

\noindent We turn next to the case $q_{0}\le N/2$, which is a consequence of the following theorem.
\begin{thm}\label{th5}
Assume $f\in C^2(\R^+)$ is nondecreasing, convex , $f>0$ in $\R^*_{+}$, \eqref{equation q} holds and $q_{0}<+\infty$.
Then, the differential inequality
\begin{equation} \label{entire subsolution} 
-\Delta u \le f(u)\qquad\text{in $\R^N$}
\end{equation}
does not admit any solution $u\in C^2(\R^N)\cap L^\infty(\R^N)$, $u>0$ such that \eqref{stability} holds, if either of the following conditions holds
\begin{enumerate}
\item $1\le N\le 9$,
\item $N=10$ and $p_{0}<+\infty$, 
\item $N\ge11$ and $p_{0}<p_{c}(N)$, 
\end{enumerate}
\end{thm}
\begin{rem}\label{lipschitz 2} 
With no change to the proof, Theorem \ref{th5} remains true  if $u$ is only assumed to be locally Lipschitz continuous. The differential inequality \eqref{entire subsolution} must then be understood in the weak sense i.e. 
$$
\int_{\R^N}\nabla u\nabla\varphi\;dx \le \int_{\R^N}f(u)\varphi\;dx,
$$
for all Lipschitz functions $\varphi\ge0$ with compact support.   
\end{rem}
It remains to prove Theorem \ref{th5}.  We begin with the following weighted-Poincar\'e inequality.
\begin{lem}\label{lemma1}
Assume $\Omega$ is an arbitrary open set in $\R^N$. Let $u\in C^2(\Omega)$, $u\ge0$ satisfy
$$
-\Delta u\le f(u)\qquad\text{in $\Omega$.}
$$
Assume in addition that for all $\varphi\in C^1_{c}(\Omega)$,
\begin{equation}
\label{stability-again}
\int_{\Omega}f'(u)\varphi^2\; dx \le \int_{\Omega} |\nabla \varphi |^2\;dx.
\end{equation}
Let $\phi\in W^{1,\infty}_{\text{\rm loc}}(\R;\R)$ denote a convex function and $\eta\in C^1_c(\R^N)$. Let 
$$\psi(u)=\int_{0}^{u}\phi'^2(t)\; dt.$$ Then,
\begin{equation}\label{crandall rabinowitz}
\int_{\Omega} [(f'\phi^2-f\psi)\circ u] \eta^2 \;dx \le \int_{\Omega} [\phi^2\circ u]\left|\nabla \eta \right|^2. 
\end{equation}
\end{lem}

\begin{rem}\label{variant}
If $\phi$ is not convex, then the following variant of \eqref{crandall rabinowitz} holds.
\begin{equation}\label{crandall rabinowitz 2}
\int_{\Omega} [(f'\phi^2-f\psi)\circ u] \eta^2 \;dx \le \int_{\Omega} [K\circ u]\Delta(\eta^2) \;dx
-\int_{\Omega} [\phi^2\circ u] \eta\Delta\eta \;dx,
\end{equation}
where $K(u)=\int_{0}^u \psi(s)\;ds$.
\end{rem}


\proof{} Multiply \eqref{entire subsolution} by $\psi(u)\eta^2$ and integrate by parts :
\begin{align*}
\int_{\Omega} \nabla u \nabla \left(\psi(u)\eta^2\right) \; dx &\le \int_{\Omega} f(u)\psi(u)\eta^2 \; dx\\
\int_{\Omega} \phi'(u)^2 |\nabla u|^2 \eta^2 \; dx + \int_{\Omega} \psi(u)\nabla u \nabla \eta^2 \; dx &\le\\
\int_{\Omega} \phi'(u)^2 |\nabla u|^2 \eta^2 \; dx - \int_{\Omega} K(u)\Delta \eta^2 \; dx &\le\\
\end{align*}
where $K(u)=\int_{0}^{u}\psi(s)\; ds$. Hence,
\begin{align}\label{integrating with psi eta squared}
\int_{\Omega} \phi'(u)^2 |\nabla u|^2 \eta^2 \; dx \le \int_{\Omega} K(u)\Delta \eta^2 \; dx + \int_{\Omega} f(u)\psi(u)\eta^2 \; dx
\end{align}
Next, we apply \eqref{stability-again} with $\varphi = \phi(u)\eta$ and obtain
\begin{align*}
\int_{\Omega} f'(u)\phi(u)^2\eta^2 \; dx &\le \int_{\Omega} |\nabla \left(\phi(u)\eta\right)|^2 \; dx = \int_{\Omega} \left|\phi'(u)\eta\nabla u + \phi(u)\nabla \eta \right|^2\; dx \\
&\le \int_{\Omega} \phi'(u)^2\eta^2|\nabla u|^2 \; dx + \int_{\Omega} \phi(u)^2|\nabla \eta|^2\; dx + 2\int_{\Omega} \phi(u)\phi'(u)\eta\nabla\eta\nabla u\; dx\\
&\le \int_{\Omega} \phi'(u)^2\eta^2|\nabla u|^2 \; dx + \int_{\Omega}\phi(u)^2|\nabla \eta|^2\; dx + \frac12\int_{\Omega} \nabla\eta^2\nabla \phi(u)^2\; dx\\
&\le \int_{\Omega} \phi'(u)^2\eta^2|\nabla u|^2 \; dx + \int_{\Omega} \phi(u)^2 \left( |\nabla \eta|^2\; -\frac12 \Delta\eta^2\right)dx \\
\end{align*}
Plug \eqref{integrating with psi eta squared} in the above. Then,
\begin{align*}
\int_{\Omega} \left( f'(u)\phi(u)^2-f(u)\psi(u)^2\right)\eta^2\; dx\le \int_{\Omega} K(u)\Delta \eta^2 \; dx +  \int_{\Omega} \phi(u)^2 \left( |\nabla \eta|^2\; -\frac12 \Delta\eta^2\right)dx
\end{align*}
This proves Remark \ref{variant}. Finally, when $\phi$ is convex,
$$
\psi(u) = \int_{0}^{u}\phi'^2(s)\;ds \le \phi'(u)\phi(u).
$$
Integrating, we obtain that $K\le \frac{1}{2}\phi^2$ and \eqref{crandall rabinowitz}  follows.
\qed

\proof{ of Theorem \ref{th5} continued}
Take $\alpha\ge1$ and $\phi=f^\alpha$. In order to take advantage of Lemma \ref{lemma1},  we need to make sure that the quantity $(f'\phi^2-f\psi)\circ u$ remains nonnegative and better, bounded below by some positive function of $u$. Clearly, the best one can hope for is an inequality of the form 
$$(f'\phi^2-f\psi)\circ u \ge c\; f'\phi^2\circ u.$$ 
To obtain such an inequality, we apply L'H\^opital's Rule :
\begin{align*}
 \lim_{0^+} \frac{f'\phi^2}{f\psi} &= \lim_{0^+} \frac{f'f^{2\alpha-1}}{\psi}\\
&= \lim_{0^+} \frac{f''f^{2\alpha-1} + (2\alpha -1)f^{2\alpha-2}f'^2}{\alpha^2 f^{2\alpha-2}f'^2}\\
&= \frac{1}{\alpha^2}\left(1/q_{0}+2\alpha-1\right)>1,
\end{align*}
where the last inequality holds if $\alpha\in [1,1+1/\sqrt{q_0})$. Note that this interval is nonempty since we assumed $q_{0}<+\infty$. Hence, for some constant $c>0$, 
\begin{equation}\label{f prime phi squared epsilon}  
f'\phi^2 - f\psi \ge c \;f'\phi^2
\end{equation}
in a neighbourhood $[0,\epsilon]$ of the origin. Modifying $\phi$, the above inequality can be extended to a given compact interval $[0,M]$ as follows. Take $\phi\in  W^{1,\infty}_{\text{\rm loc}}(\R;\R)$ defined by
\begin{equation}\label{newphi} 
\phi(u)=\left\{ 
\begin{aligned} 
&f(u)^\alpha&\text{ if $0\le u\le \epsilon $}          \\ 
&f(\eps)^{\alpha-1}f(u)\exp\left(\int_{\eps}^{u}\sqrt{\frac{f''}{f}}\;ds\right)&\text{ if $u>\epsilon$}          \\ 
\end{aligned} 
\right.
\end{equation}
where  $\epsilon, \alpha$ are chosen as before. Then $\phi\in W^{1,\infty}_{\text{\rm loc}}(\R;\R)$.  For $u>\eps$, we claim that the quantity $\frac{f'}{f}\phi^2 - \psi$ is constant. Indeed,
\begin{align*}
\left(\frac{f'}{f}\phi^2-\psi\right)' &=  \left(\frac{f'}{f}\right)'\phi^2 + 2\frac{f'}{f}\phi\phi' - \phi'^2\\
&= \left(\frac{f''}{f}-\frac{f'^2}{f^2}\right)\phi^2 + 2\frac{f'}{f}\phi\phi' - \phi'^2\\
&= \frac{f''}{f}\phi^2 - \left(\frac{f'}{f}\phi - \phi'\right)^2 = \phi^2\left(\frac{f''}{f}-\left(\frac{f'}{f} - \frac{\phi'}{\phi}\right)^2\right)\\  
&=0, 
\end{align*}
where we used the definition of $\phi$ in the last equality. So for $u>\eps$, 
\begin{align*} 
f'\phi^2 - f\psi &= f\left(\frac{f'}{f}\phi^2-\psi\right) = f\left(\frac{f'(\eps)}{f(\eps)}\phi^2(\eps)-\psi(\eps)\right)\\
&\ge f'(\eps)\phi^2(\eps)-f(\eps)\psi(\eps)
\ge c_{\eps}>0,
\end{align*} 
where we used \eqref{f prime phi squared epsilon} at $u=\eps$. Since $f'\phi^2$ is bounded above by a constant on any compact interval of the form $[\eps,M]$, we conclude that \eqref{f prime phi squared epsilon} holds throughout $[0,M]$ for a constant $c>0$ perhaps smaller. 
We have just proved that given $\alpha\in [1,1+1/\sqrt{q_0})$ and a bounded positive function $u$, there exists $c>0$ such that
\begin{equation}\label{f prime phi squared old} 
[f'\phi^2 - f\psi]\circ u \ge c \;[f'\phi^2]\circ u.
\end{equation}
Recall that we established the above inequality in order to apply Lemma \ref{lemma1}. Unfortunately, since the function $\phi$ we introduced in \eqref{newphi} may not be convex, we cannot apply Lemma \ref{lemma1} directly. We make use of \eqref{crandall rabinowitz 2} instead. In order to obtain a meaningful result, we need to understand how the different functions of $u$ introduced in \eqref{crandall rabinowitz 2} compare. By definition of $\phi$, we easily deduce the following set of inequalities
\begin{equation}\label{f prime phi squared} 
\left\{ 
\begin{aligned}
\;[ f'\phi^2 - f\psi ] \circ u\ge c\, f'f^{2\alpha}\circ u\\
\;\phi^2\circ u \le C f^{2\alpha}\circ u\\
\;K\circ u \le C f^{2\alpha}\circ u,
\end{aligned}
\right.
\end{equation} 
So, we just need to relate $f$ and $f'$ to be able to compare all quantities involved in the estimate. Fix $q_{1}<q_{0}$. By definition of $q_{0}$, there exists a neighborhood of zero where
\begin{equation*}
 \frac{ff''}{f'^2}\le 1/q_{1}.
\end{equation*}
In particular, $f'/f^{1/q_{1}}$ is nonincreasing and in a neighborhood of zero we have
\begin{equation}\label{f prime bounded below}
 f'\ge cf^{1/q_{1}}.
\end{equation}
By continuity, up to choosing $c>0$ smaller, the above inequality holds in the whole range of a given bounded positive function $u$. 
Recall now \eqref{f prime phi squared}, \eqref{f prime bounded below}  and apply \eqref{crandall rabinowitz 2}. The estimate reduces to
\begin{equation*} 
 \int_{\R^N} [f^{1/{q_{1}}+2\alpha}\circ u] \eta^{2} \; dx\le C \int_{\R^N}  [f^{2\alpha}\circ u] \left(|\nabla\eta|^2+ \left| \eta\Delta\eta\right|\right)  \; dx
\end{equation*} 
Choose $\eta=\zeta^m$, $m\ge1$, $\zeta\in C^2_c(\R^N)$, $1\ge\zeta\ge0$ :
\begin{align*}
 \int_{\R^N} [f^{1/{q_{1}}+2\alpha}\circ u] \zeta^{2m} \; dx&\le C \int_{\R^N}  [f^{2\alpha}\circ u] \left(\zeta^{2m-2}|\nabla\zeta|^2+ \zeta^{2m-1}\left| \Delta \zeta\right| \right)\; dx\\
&\le C \int_{\R^N}  [f^{2\alpha}\circ u] \zeta^{2m-2}\left(|\nabla\zeta|^2+ \left| \Delta \zeta\right| \right)\; dx
\end{align*}
Using  H\"older's inequality, it follows that
\begin{align*} 
\int_{\R^N}  [f^{1/q_{1}+2\alpha}\circ u]\zeta^{2m}\; dx&\le C\left(\int_{\R^N}  [f^{2\alpha m'}\circ u]\zeta^{2m}\; dx\right)^{1/m'} \left(\int_{\R^N}\left(|\nabla\zeta|^{2}+\left| \Delta\zeta \right| \right)^m\;\right)^{1/m}.\\
\end{align*}
Assume temporarily that
\begin{equation}\label{test alpha m prime}
 f^{1/q_{1}+2\alpha}\circ u \ge c \; f^{2\alpha m'}\circ u.
\end{equation}
Then, the inequality simplifies to
\begin{align*} 
\int_{\R^N}  [f^{2\alpha m'}\circ u]\zeta^{2m}\; dx&\le C\int_{\R^N}\left(|\nabla\zeta|^{2}+\left| \Delta\zeta \right| \right)^m.
\end{align*}
Choose now $\zeta$ such that $\zeta\equiv 1$ in $B_R$ and $|\nabla\zeta|\le C/R$, $\left| \Delta\zeta\right|\le C/R^2 $ :
\begin{equation}\label{decay}
 \int_{\R^N}  [f^{2\alpha m'}\circ u] \zeta^{2m} \;dx\le C R^{N-2m},
\end{equation}
The above inequality is true as soon as \eqref{test alpha m prime} holds, which itself reduces to choosing the exponents such that
$$2\alpha m'\ge 1/q_{1}+2\alpha.$$ 
This holds for some $q_{1}<q_{0}$ provided $2\alpha(m'-1)>1/q_{0}$. Since $\alpha$ can be chosen arbitrarly close to $1+1/\sqrt{q_{0}}$  and restricting to $m'$ less than but as close as we wish to $\frac{N}{N-2}$, we finally need only assume
\begin{equation}\label{final algebra}
 \frac{4}{N-2}\left(1+1/\sqrt{q_{0}}\right)>1/{q_{0}}.
\end{equation}
Since $m'<\frac{N}{N-2}$, $N-2m<0$. So, the right-hand-side of \eqref{decay} converges to $0$ as $R\to\infty$, whence $f\circ u=0$ and $u=0$, as desired. Solving \eqref{final algebra} for $q_{0}$ yields the conditions stated in Theorem \ref{th5}. 
\hfill\qed

\section{Extensions to unbounded and sign-changing solutions}
We deal first with possibly unbounded solutions.
\proof{ of Corollary \ref{unbounded}}
Note that by Lemma \ref{lemma2}, we need only consider the case $q_{0}<+\infty$. We modify the rest of the proof of Theorem \ref{th5}  as follows : take $\phi\in  W^{1,\infty}_{\text{\rm loc}}(\R;\R)$ defined by
\begin{equation*}
\phi(u)=\left\{ 
\begin{aligned} 
&f(u)^\alpha&\text{ if $0\le u\le \epsilon $}          \\ 
&f(\eps)^{\alpha-1}f(u)\exp\left(\int_{\eps}^{u}\sqrt{\frac{f''}{f}}\;ds\right)&\text{ if $\epsilon < u\le 1/\epsilon$}          \\ 
&f(u)^\beta+A&\text{ if $u> 1/\epsilon$}, 
\end{aligned} 
\right.
\end{equation*}
where $\alpha$ is chosen in $[1,1+1/\sqrt{q_0})$ as previously, $\beta$ in $[1,1+1/\sqrt{\overline{q_{\infty}}})$ and $A$ such that $\phi$ is $W^{1,\infty}_{\text{\rm loc}}(\R;\R)$.  
Then, \eqref{f prime phi squared old} holds if in addition 
$$\liminf_{u\to+\infty}\frac{f'\phi^2}{f\psi}(u)>1.$$ 
We leave the reader check that this is true under assumption \eqref{equation q infty}, for $\beta\in[1,1+1/\sqrt{\overline{q_{\infty}}})$. Apply \eqref{crandall rabinowitz 2}  with $\eta = \zeta^m$, $m\ge 1$, $\zeta\in C^2_c(\R^N)$, $0\le\zeta\le1$ :
\begin{multline}\label{phiphi} 
 \int_{\R^N} [f'\phi^2\circ u] \zeta^{2m} \; dx \le C \int_{\R^N}  [\left(\phi^{2}+K\right)\circ u] \zeta^{2m-2}\left(|\nabla\zeta|^2+ \left| \Delta\zeta \right|\right) \; dx
\\
 \le C \left( \int_{\R^N}[\left(\phi^{2}+K\right)^{m'} \circ u] \; \zeta^{2m}\; dx \right)^{1/m'} 
\left( \int_{\R^N} \left(\left| \nabla \zeta\right|^2 + \left| \Delta\zeta \right| \right)^m \; \; dx \right)^{1/m}
\end{multline}
By definition of $\phi$ and \eqref{f prime bounded below}, there exists constants $c,c'>0$ such that for $u\in [0,1]$, $f'\phi^2 (u) \ge c f'f^{2\alpha}(u) \ge c' f^{2\alpha +q_{1}}(u)$, where $q_{1}<q_{0}$. We also clearly have $\left(\phi^{2}+K\right)^{m'}(u)\le C f^{2\alpha m'}$ for $u\in [0,1]$. So, 
$$f'\phi^2\ge c\left(\phi^{2}+K\right)^{m'}\qquad\text{ on $[0,1]$},$$ provided that $2\alpha m'\ge 1/q_{1}+2\alpha$. Similarly, the reader will easily check using \eqref{equation q infty} that given $q_{2}<\overline{q_{\infty}}$, there exists $c>0$ such that
$$
f'\ge c f^{1/q_{2}}\qquad\text{in $[1,+\infty)$},
$$ 
whence $f'\phi^2\ge c\left(\phi^{2}+K\right)^{m'}$ in $[1,+\infty)$ provided that $2\beta m'\ge 1/q_{2}+2\alpha$. We conclude that
\begin{equation} \label{summary}
f'\phi^2\circ u\ge c\left(\phi^{2}+K\right)^{m'}\circ u,
\end{equation}  
provided that $2\alpha m'\ge 1/q_{1}+2\alpha$ and $2\beta m'\ge 1/q_{2}+2\alpha$.
Since $\alpha$ can be chosen arbitrarily close to $1+1/\sqrt{q_{0}}$, $\beta$ to $1+1/\sqrt{\overline{q_{\infty}}}$, $q_{1}$ to $q_{0}$, $q_{2}$ to $\overline{q_{\infty}}$ and $m'$ to $N/(N-2)$, we conclude that suitable parameters can be chosen provided \eqref{final algebra} and 
$$
 \frac{4}{N-2}\left(1+1/\sqrt{\overline{q_{\infty}}}\right)>1/\overline{q_{\infty}}
$$
 hold. These inequalities are true under the assumptions of Remark \ref{unbounded}. So, collecting \eqref{phiphi} and \eqref{summary}, we obtain for some $m>N/2$,
\begin{align}\label{phiphiphi} 
 \int_{\R^N} [\left(\phi^{2}+K\right)^{m'}\circ u] \zeta^{2m} \; dx&\le C\int_{\R^N} \left(\left| \nabla \zeta\right|^{2}+ \left| \Delta\zeta \right| \right)^{m'} \; dx.
\end{align}  
 Choose at last $\zeta$ such that $\zeta\equiv 1$ in $B_R$ and $|\nabla\zeta|\le C/R$, $\left| \Delta\zeta\right|\le C/R^2 $: the right-hand side of \eqref{phiphiphi} converges to $0$ as $R\to\infty$ and the conclusion follows. \qed
We work next with sign-changing solutions.
\proof{ of Corollaries \ref{sign changing} and \ref{sign changing 2}}
We simply remark that if $u\in C^2(\R^N)$ is a solution of \eqref{main}, then $u^+$ (respectively $u^-$) is locally Lipschitz continuous and solves the differential inequality \eqref{entire subsolution} (respectively $-\Delta u^-\le \tilde f(u^-)$ in $\R^N$, where $\tilde f(t):=-f(-t)$ for $t\in\R^-$). Since we assumed $q_{0}<+\infty$, we may then apply Theorem \ref{th5} and Remark \ref{lipschitz 2}.  Corollary \ref{sign changing} follows. For Corollary \ref{sign changing 2}, we replace Theorem \ref{th5} by the adaptation presented in the proof of Corollary \ref{unbounded}.   
\hfill\qed

\section{Beyond power-type nonlinearities}
\proof{ of Theorem \ref{thbeyond}}
Case \ref{case 1}. of the theorem was proved in Lemma \ref{lemma2}. For cases \ref{case 2}. and \ref{case 3}.       
take $\alpha\ge1$ and $\phi=f^\alpha$. Let $L=\liminf_{0^+}{f'\phi^2}/{f\psi}$ and let $(u_{n})$ denote a sequence along which ${f'\phi^2}/{f\psi}$ converges to $L$. 
By Remark \ref{remark pohozaev},  we may always assume that $f(0)=0$. So, applying Cauchy's mean value theorem, there exists $v_{n}\in(0,u_{n})$ such that
\begin{align*}
\frac{f'\phi^2}{f\psi}(u_{n})  &= \frac{f'f^{2\alpha-1}}{\psi}(u_{n})\\
&=\left.\frac{f''f^{2\alpha-1} + (2\alpha -1)f^{2\alpha-2}f'^2}{\alpha^2 f^{2\alpha-2}f'^2}\right|_{u=v_{n}}
\end{align*}
Passing to the limit, we obtain
\begin{equation}\label{theliminf}  
L= \liminf_{0^+} \frac{f'\phi^2}{f\psi}\ge\frac{1}{\alpha^2}\left(\frac1{\overline{q_{0}}}+2\alpha-1\right)>1,
\end{equation} 
where the last inequality holds if $\alpha\in [1,1+1/\sqrt{\overline{q_{0}}})$. Note that this interval is nonempty since we assumed ${\overline{q_{0}}<\infty}$. At this point, we repeat the steps performed in the proof of Theorem \ref{th5} : from equation \eqref{theliminf}, we deduce that \eqref{f prime phi squared epsilon} holds in a neighborhood $[0,\eps]$ of the origin. Modifying $\phi$ as in \eqref{newphi}, the \it verbatim \rm arguments lead to \eqref{f prime phi squared old} and \eqref{f prime phi squared}. For the rest of the proof, we argue slightly differently according to the case considered.

\noindent {\bf Case \ref{case 2}. of Theorem \ref{thbeyond}}    
In place of \eqref{f prime bounded below}, we simply use the convexity of $f$. Since $u$ is bounded, there exists a constant $c>0$ such that
$$
f'(u)\ge\frac{f(u)}{u}\ge c f(u).
$$  
So, \eqref{decay} holds for some $m>N/2$ whenever $\frac{4}{N-2}\left(1+1/\sqrt{\overline{q_{0}}}\right)>1$, which is true for $N\le 6$, provided $\overline{q_{0}}<\infty$. Following the proof of Theorem        
\ref{th5}, we obtain case \ref{case 2}.  of Theorem \ref{thbeyond}.     

\noindent {\bf Case \ref{case 3}. of Theorem \ref{thbeyond}} 
By definition of $\underline{q_{0}}$, \eqref{f prime bounded below} now holds for $q_1<\underline{q_{0}}$. Resuming our inspection of the proof of Theorem        
\ref{th5}, we see that \eqref{decay} holds under assumption \ref{case 3}. of Theorem \ref{thbeyond} and the desired conclusion follows.     
\hfill\qed
 
\section{Speed of decay : proof of Theorem \ref{th3}  }
In this section, we characterize the speed of decay of solutions which are stable outside a compact set. To do so, we shall again take advantage of  Lemma \ref{lemma1} or actually its general form \eqref{crandall rabinowitz 2}, with a different choice of test function $\phi\circ u$. We divide the proof in several steps.

\noindent{\bf Step 1. } We begin by proving the usual estimate
\begin{equation*}
[f'\phi^2 - f\psi] (u) \ge c [f'\phi^2](u)
\end{equation*} 
where this time $\phi(u)=\left(\frac{f(u)}{u}\right)^\alpha$ and $\alpha$ is chosen in a suitable range.

First, by Lemma \ref{lemma2} and Remark \ref{remark pohozaev} , we may restrict to the case where $q_{0}<+\infty$, whence  $p_{0}>1$, and we may also assume $f(0)=f'(0)=0$. By Proposition \ref{prop1}, $\lim_{\left|  x \right|\to\infty }u(x)=0$.  For $u\in\R^*_{+}$, take $\phi\in  W^{1,\infty}_{\text{\rm loc}}(\R;\R)$ defined by
\begin{equation}\label{new def of phi} 
\phi(u)=\left(\frac{f(u)}{u}\right)^\alpha,
\end{equation}
where $\alpha>\frac12$. 
We begin by computing
\begin{equation} \label{quotient} 
L = \liminf_{u\to0^+}\frac{f'\phi^2}{f\psi} (u).
\end{equation} 
Let $(u_{n})$ denote a sequence along which ${f'\phi^2}/{f\psi}$ converges to $L$. 
Observe that since $f(0)=f'(0)=0$, then $\psi(0)=0$ and 
\begin{equation}\label{alpha bigger than one half} 
 \lim_{u\to0}f'f^{2\alpha-1}u^{-2\alpha}= \lim_{u\to0}\frac{f'(u)}{u}\left(\frac{f(u)}{u}\right)^{2\alpha-1}=0,
\end{equation}  
if $\alpha>1/2$. So, 
applying Cauchy's mean value theorem, there exists $v_{n}\in(0,u_{n})$ such that
\begin{align*}
\frac{f'\phi^2}{f\psi}(u_{n})  &= \left.\frac{f'f^{2\alpha-1}u^{-2\alpha}}{\psi} 
\right|_{u=u_{n}}\\
&=\left. \frac{f''f^{2\alpha-1}u^{-2\alpha}  + (2\alpha-1)f'^2f^{2\alpha-2}u^{-2\alpha} - 2\alpha f'f^{2\alpha-1}u^{-2\alpha-1} }{\alpha^2u^{-2\alpha-2}f^{2\alpha}(-1+uf'/f)^2} \right|_{u=v_{n}} \\
&=  \left. \frac{{f''}u^2/{f}\,  + (2\alpha-1)f'^2u^2/f^2 -2\alpha uf'/f}{\alpha^2(-1+uf'/f)^2}  \right|_{u=v_{n}}\\
&= \left. \frac{ff''/f'^2 + (2\alpha-1)-2\alpha f/(uf') }{\alpha^2 (1-f/(uf'))^2}\right|_{u=v_{n}}
\end{align*}
For $u\in\R^*_{+}$, let
\begin{equation}\label{equation p}
p(u) = \frac{uf'(u)}{f(u)}
\end{equation}
It follows that
\begin{align}\label{quotient simplified} 
\frac{f'\phi^2}{f\psi}(u_{n})   &= \left. \frac{1/q+2\alpha-1-2\alpha/p}{\alpha^2(1-1/p)^2}  \right|_{u=v_{n}}  \\
&=  1 + \left. \frac{1/q-(\alpha(1-1/p) - 1)^2}{\alpha^2(1-1/p)^2}  \right|_{u=v_{n}} \nonumber
\end{align}
We claim that \eqref{equation q} implies
\begin{equation}
\label{p goes to pnought}
p_{0} = \lim_{u\to 0^+} p(u),
\end{equation}
where $p_{0}$ is the conjugate exponent of $q_{0}$ i.e. \eqref{pnought} holds.
Take indeed any cluster point $p_{1}$ of $p$ and a sequence $(u_{n})$ such that $p$ converges to $p_{1}$ along $(u_{n})$. Applying Cauchy's mean value theorem, there exists $v_{n}\in (0,u_{n})$ such that
$$
p(u_{n}) = \left. \frac{f'+uf''}{f'}\right|_{u=v_{n}} = 1 + p/q (v_{n}).  
$$
Let $\underline{p_{0}} = \liminf_{u\to 0^+} p(u)$ and $\overline{p_{0}} = \limsup_{u\to 0^+} p(u)$. Pass to the limit as $n\to+\infty$ :
$$
1+\subp/q_{0} \le p_{1} \le 1+ \superp/q_{0}.
$$
Applying the above inequality to $p_{1}=\subp, \;\superp$, we obtain
$$
\superp(1-1/q_{0})\le1\le \subp(1-1/q_{0})
$$
and \eqref{p goes to pnought} follows. Next, we apply \eqref{p goes to pnought} in \eqref{quotient simplified}. Thus,
$$
L = 1 + \frac{1/q_{0} - (\alpha/q_{0} -1)^2}{\alpha^2/q_{0}^2}.
$$
So, $L>1$ if 
\begin{equation} \label{range alpha} \alpha\in (q_{0}-\sqrt {q_{0}},q_{0} + \sqrt{q_{0}}).\end{equation} 
We conclude that given $\alpha> 1/2$ in the range \eqref{range alpha}, there exists $c>0$ such that for $u$ small enough 
\begin{equation} \label{f of u over u} 
[f'\phi^2 - f\psi] (u) \ge c [f'\phi^2](u)\ge c \left(\frac{f(u)}{u}\right)^{2\alpha+1},
\end{equation} 
where we used the convexity of $f$ in the last inequality. Note that since $u(x)\to0$ as $\left| x \right| \to+\infty$, the above inequality holds for $u=u(x)$ and $x$ in the complement of a ball of large radius.

\noindent{\bf Step 2. }
Next, we need to estimate the other functions of $u$ appearing in \eqref{crandall rabinowitz 2}. We claim that for small values of $u$,
\begin{equation} \label{estimate K}
K(u)\le C \left(\frac{f(u)}{u}\right)^{2\alpha}.
\end{equation} 
To see this, it suffices to prove that $\limsup_{u\to0^+}K(u)/\phi^2(u)<\infty$. Take a sequence $(u_{n})$ converging to zero and apply Cauchy's mean value theorem : there exists $v_{n}\in (0,u_{n})$ such that
\begin{align*}
\frac{K}{\phi^2}(u_{n})  &= \frac{\psi}{2\phi\phi'}(v_{n}) 
\end{align*}
It follows from \eqref{f of u over u} that $f'\phi^2-f\psi\ge0$ for small $u$. So, $\psi(v_{n})\le [f'\phi^2/f](v_{n})$ for large $n$ so that
\begin{align*}
\frac{K}{\phi^2}(u_{n})  &\le \frac{f'\phi}{2f\phi'}(v_{n})=\frac{1}{2\alpha(1-1/p(v_{n}))}. 
\end{align*}
Recalling \eqref{p goes to pnought} and since we assumed that $p_{0}>1$, \eqref{estimate K} follows.  

\noindent{\bf Step 3. }
In this step, we prove an estimate of the form
\begin{equation*}
 \int_{B_{R}(x_{0})} \left(\frac{f(u)}{u}\right)^{2\alpha+1} \; dx
 \le C R^{N-2m}, 
\end{equation*} 
where $m=2\alpha+1$ and $B_{R}(x_{0})$ is a suitably chosen ball shifted towards infinity. 
Choose  $\zeta\in C^2_{c}(\R^N)$, $0\le\zeta\le1$ supported outside a ball $B_{R_{0}}(0)$ of large radius, so that \eqref{stability} holds for functions supported outside $B_{R_{0}}(0)$ and that 
\eqref{f of u over u} and \eqref{estimate K} hold for $u=u(x)$, $x\in \text{\rm supp}\,\zeta$. 
By Lemma \ref{lemma1}, we may apply \eqref{crandall rabinowitz 2} with $\eta=\zeta^m$, $m\ge1$. Using \eqref{f of u over u}, \eqref{estimate K} and the convexity of $f$
, we obtain for $\alpha>1/2$ in the range \eqref{range alpha} 
\begin{align*} 
 \int_{\R^N} \left(\frac{f(u)}{u}\right)^{2\alpha+1} \zeta^{2m}\; dx&\le \int_{\R^N} f'(u) \left(\frac{f(u)}{u}\right)^{2\alpha} \zeta^{2m} \; dx\\
 &\le C \int_{\R^N}   \left(\frac{f(u)}{u}\right)^{2\alpha}\zeta^{2m-2} \left(|\nabla\zeta|^2+ \left| \Delta\zeta\right|\right)  \; dx.\\
\end{align*}
Fix $m=2\alpha+1$ and apply H\"older's inequality. It follows that
\begin{align}\label{ineq zeta}  
 \int_{\R^N} \left(\frac{f(u)}{u}\right)^{2\alpha+1} \zeta^{2m}\; dx
 &\le C \int_{\R^N} \left(|\nabla\zeta|^2+ \left| \Delta\zeta\right|\right)^{2m}  \; dx.
\end{align}      
We work on balls shifted towards infinity. More precisely, we take a point $x_{0}\in\R^N$ such that $\left| x_{0}\right|>10R_{0}$ and set $R= \left| x_{0} \right|/4$. Then, $B(x_0,2R)\subset \{ x\in \R^N\; : \; \left| x \right| \ge R_{0}\}$ and
we may apply \eqref{ineq zeta} with $\zeta=\varphi(\left| x-x_{0} \right|/R)$ and $\varphi\in C^2_{c}(\R)$ given by
\begin{equation*}
\varphi(t)=\left\{ 
\begin{aligned}
1&\qquad\text{if $\left| t\right|\le1$, }\\
0&\qquad\text{if $\left| t\right|\ge 2$.}
\end{aligned}
\right.
\end{equation*} 
We get
\begin{equation}\label{second estimate} 
 \int_{B_{R}(x_{0})} \left(\frac{f(u)}{u}\right)^{2\alpha+1} \; dx
 \le C_{3}R^{N-2m}. 
\end{equation} 
\noindent{\bf Step 4. } In this step, we prove the estimate
$$
R^\epsilon \| f(u)/u \|_{L^{\frac N{2-\eps}}(B(x_{0},R))} \le C
$$
By Lemma \ref{basic lemma}, $q_{0}\ge1$. Under the assumptions of Theorem \ref{th3}, we can choose the exponent $m$ so large that for small $\eps>0$, $m>N/(2-\eps)$ (recall that $m=2\alpha+1$ and $\alpha>1/2$ can be chosen freely in the range \eqref{range alpha}) . 
Furthermore, by H\"older's inequality and \eqref{second estimate}, we obtain     
\begin{align} \label{control serrin} 
R^\epsilon \| f(u)/u \|_{L^{\frac N{2-\eps}}(B(x_{0},R))} &\le R^\epsilon \| f(u)/u \|_{L^{m}(B(x_{0},R))}\left| B_{R}\right|^{\frac{2-\eps}{N}-\frac1m}   \\ 
&\le C R^{\eps}\left(R^{N-2m}\right)^{1/m}R^{2-\eps-\frac Nm}= C.\nonumber
\end{align} 

\noindent{\bf Step 5. }
Now, we think of $u$ as a solution of a linear problem, namely
\begin{equation} \label{linear serrin} 
-\Delta u = \frac{f(u)}{u}\, u=:V(x)u\qquad\text{in $\R^N$.}
\end{equation} 
According to classical results of J. Serrin \cite{serrin} and N. Trudinger \cite{trudinger} (see also Theorem 7.1.1 on page 154 of \cite{pucci-serrin}), for any $p\in(1,+\infty)$ and any $x_{0}\in\R^N$, there exists a constant 
\begin{equation*} 
C_{S}=C_{S}(R^\eps \| V\|_{L^{\frac N{2 -\eps}}(B(x_{0},2R))}, N,p)>0
\end{equation*}  
such that 
\begin{equation} \label{serrin inequality} 
\|u \|_{L^\infty(B_{R}(x_{0}))} \le C_{S} R^{-N/p} \| u\|_{L^p(B_{2R}(x_{0}))}.  
\end{equation}
Note that for our choice of $x_{0}$, equation \eqref{control serrin} holds and so $C_{S}$ is a true constant, independant of $R$ and $x_{0}$.

\noindent{\bf Step 6. } 
The inequality \eqref{serrin inequality} gives a pointwise estimate in terms of an integral average of $u$. In order to control the latter, we consider $\tilde u$ the average of $u$ over the sphere $\partial B_r(x_{0})$, defined for $r>0$ by $\tilde u(r) = \fint_{\partial B_{r}(x_{0})}u\;d\sigma$. We claim that there exists $C=C(N)>0$ such that
\begin{equation}\label{estimate u tilde}
\frac{f(\tilde u(r))}{\tilde u(r)}\le \frac C{r^2}.
\end{equation}
To prove this, we first observe that since $f$ is convex, $\tilde u$ satisfies the differential inequality
$$
- \tilde u'' - \frac{N-1}{r}\tilde u' \ge f(\tilde u).
$$
Now,  since $f\ge0$, $\tilde u'\le0$. In particular $r\mapsto f(\tilde u(r))$ is nonincreasing. Fix $\lambda\in (0,1)$ and integrate the differential inequality between $0$ and $r$ :
\begin{equation*} 
-\tilde u'(r) \ge r^{1-N}\int_{0}^r s^{N-1}f(\tilde u(s))\;ds\ge r^{1-N}\int_{0}^{\lambda r}s^{N-1}f(\tilde u(s))\;ds \ge \frac{\lambda^Nrf(\tilde u(\lambda r))}{N}.
\end{equation*}   
Integrate a second time between $r$ and $r/\lambda$. Then,
\begin{equation*}
\tilde u(r) \ge \tilde u(r/\lambda) + \frac{\lambda^N}{N}\int_{r}^{r/\lambda}sf(\tilde u(\lambda s))\;ds\ge r^2f(\tilde u(r)) \frac{\lambda^N}{2N}\left(\frac{1}{\lambda^2} -1\right).
\end{equation*}
Taking $\lambda=\frac{N-2}{N}$, \eqref{estimate u tilde}  follows with $C=\frac{1}{2N}\left(\frac{N-2}{N}\right)^N\left(\left(\frac{N}{N-{2}}\right)^{2} -1\right)$.  

\noindent{\bf Step 7. }Recall that we are trying to establish an $L^p$ estimate, $p>1$ in order to use \eqref{serrin inequality} . To start with, we use \eqref{estimate u tilde} to obtain an $L^1$ estimate of $f(u)$. Namely, we prove that there exist constants $C_{1}, C_{2}>0$ depending on $N$ only, such that
\begin{equation}\label{estimate average f} 
\fint_{B_{R}(x_{0})} f(u)\;dx \le C_{1} R^{-2} g(C_{2}/R^2),
\end{equation}
where $g$ is the inverse function of $t\mapsto\frac{f(t)}{t}$, which exists for small values of $t$ by Remark \ref{f sur t est inversible}.  For simplicity, we write $B_{R}$ in place of $B_{R}(x_0)$ in what follows. To prove \eqref{estimate average f} , observe that for $r\in (R,2R)$,
\begin{align*}
\int_{B_{R}} f(u)\; dx = &c_{N} R^{N-2}\int_{R}^{2R}r^{1-N}\;dr \int_{B_{R}} f(u)\; dx\\
 \le &c_{N} R^{N-2}\int_{R}^{2R}r^{1-N}\;dr \int_{B_{r}} f(u)\; dx
= \\
&c_{N} R^{N-2}\int_{R}^{2R}r^{1-N}\;dr \int_{B_{r}} -\Delta u\; dx\\
\le& -c_{N} R^{N-2}\int_{R}^{2R}r^{1-N}\;dr \int_{\partial B_{r}} \frac{\partial u}{\partial n}\; d\sigma =\\
& -c_{N} R^{N-2}\int_{R}^{2R}\tilde u'\;dr\le c_{N}R^{N-2}\tilde u(R).\\
\end{align*}
Estimate \eqref{estimate average f} follows, using \eqref{estimate u tilde}. 

\noindent{\bf Step 8. }The assumptions on $f$ allow us to convert \eqref{estimate average f} into an $L^p$ estimate.  Indeed, since $q_{0}<\infty$ (in fact, one only needs $\overline{q_{0}}<\infty$), one can easily check that there exists $p>1$ such that the function $h(t)=f(t^{1/p})$ is convex for small $t$. By Jensen's inequality,
\begin{equation*} 
h\left(\fint_{B_{R}}u^p\;dx\right) \le \fint_{B_{R}}f(u)\;dx \le C_{1}R^{-2}g(C_{2}/R^2).
\end{equation*} 
By Remark \ref{f sur t est inversible}, $f$ is invertible and so is $h$. Composing by $h^{-1}$, we obtain
\begin{equation*}
\int_{B_{R}}u^p\;dx\le C R^N h^{-1}\left(C_{1}R^{-2}g(C_{2}/R^2)\right).
\end{equation*}
Combining this with \eqref{serrin inequality}, we finally obtain
\begin{multline*} 
\| u \|_{L^\infty(B_{R})}\le  C R^{-N/p}\left(R^N h^{-1}\left(C_{1}R^{-2}g(C_{2}/R^2)\right)\right)^{1/p}=\\
C f^{-1}\left(C_{1}R^{-2}g(C_{2}/R^2)\right)=C\; s(R).
\end{multline*}  
\hfill\qed


\noindent We conclude this section by proving Corollary \ref{cor fine asymptotics}. Namely, we improve the rate of decay from $O(s( \left|  x \right| ))$ to $o(s( \left| x \right| ))$, when additional information on the nonlinearity is available. 
\proof{ of Corollary \ref{cor fine asymptotics}}  
To start with, observe that under assumption \eqref{lower and upper bound}, there exists a constant $C>0$ such that
\begin{equation} \label{sharp rate of decay}
s(R)\le CR^{-\frac2{p_{0}-1}}.
\end{equation}  
Recall now \eqref{ineq zeta}. We choose a suitable cut-off function $\zeta\in C^2_{c}(\R^N)$ as follows. Let $\varphi\in C^2_{c}(\R)$ satisfying $0\le\varphi\le1$ everywhere on $\R$ and
\begin{equation*}
\varphi(t)=\left\{ 
\begin{aligned}
1&\qquad\text{if $\left| t\right|\le1$, }\\
0&\qquad\text{if $\left| t\right|\ge 2$.}
\end{aligned}
\right.
\end{equation*} 
For $s>0$, let $\theta_{s}\in C^2_{c}(\R)$ satisfying $0\le\theta_{s}\le1$ everywhere on $\R$ and
\begin{equation*}
\theta_{s}(t)=\left\{ 
\begin{aligned}
0&\qquad\text{if $\left| t\right|\le s+1$, }\\
1&\qquad\text{if $\left| t\right|\ge s+2$.}
\end{aligned}
\right.
\end{equation*} 
Given $R>R_{0}+3$, we define $\zeta$ at last by
\begin{equation*}
\zeta(x)=\left\{ 
\begin{aligned}
\theta_{R_{0}}( \left| x \right|) &\qquad\text{if $\left| x\right|\le R_{0}+3$, }\\
\varphi( \left| x \right|/R )&\qquad\text{if $\left| x\right|\ge R_{0}+3$.}
\end{aligned}
\right.
\end{equation*} 
Applying \eqref{ineq zeta} with $\zeta$ as above, we deduce that for some constants $C_1,C_{2}>0$,
\begin{equation}\label{first estimate} 
 \int_{B_{R}\setminus B_{R_0+2}} \left(\frac{f(u)}{u}\right)^{2\alpha+1} \; dx
 \le C_{1} + C_{2}R^{N-2m}. 
\end{equation}
Recall that \eqref{first estimate} holds for $m=2\alpha+1$ and any $\alpha>1/2$ such that $q_{0}-\sqrt{q_{0}}<\alpha<q_{0}+\sqrt{q_{0}}$. In fact, the restriction $\alpha>1/2$ can be lifted and replaced by $\alpha>0$. Indeed,  the restriction $\alpha>1/2$ was used for the sole purpose of proving \eqref{alpha bigger than one half}.  But \eqref{alpha bigger than one half} clearly holds under the finer assumption  \eqref{lower and upper bound} for any $\alpha>0$.

We would like to choose $\alpha$ such that $m:=2\alpha+1=N/2$. 
Since $p_{0}$ is in the supercritical range \eqref{supercritical range} , straightforward algebraic computations show that such a choice is indeed possible in the range $q_{0}-\sqrt{q_{0}}<\alpha<q_{0}+\sqrt{q_{0}}$. By \eqref{first estimate}, we deduce that
$$
\int_{\R^N}u^{(p_{0}-1)\frac N2}<\infty.
$$ 
In particular, given $\eta>0$ small, there exists $R>0$ so large that given any point $x_{0}\in\R^N$ such that $\left| x_{0} \right| =4R$,
$$
\int_{B_{R}(x_{0})}u^{(p_{0}-1)\frac N2}<\eta.
$$ 
We apply again \eqref{serrin inequality}, this time with $p=(p_{0}-1)\frac N2$ and obtain
\begin{equation}\label{u is small o} 
\|u \|_{L^\infty(B_{R}(x_{0}))} \le C_{S} R^{-N/p} \| u\|_{L^p(B_{2R}(x_{0}))}\le C_{S}\eta R^{-\frac2{p_{0}-1}}.  
\end{equation}
This shows that $u(x)=o(\left| x \right|^{-\frac2{p_{0}-1}} )$.
It remains to prove the estimate on $\left| \nabla u \right|$. 
Observe that any partial derivative $v=\partial u / \partial x_{i}$ solves the linearized equation
$$
-\Delta v = f'(u)\, v\qquad\text{in $\R^N$.}
$$
Apply again the Serrin inequality \eqref{serrin inequality}, this time with potential $\tilde V(x)=f'(u)$ and solution $v$. Since $0\le f'(u)\le C u^{{p_{0}-1}}$, the potential $\tilde V$ is equivalent to $V(x)=f(u)/u$ and so the Serrin constant $C_{S}$ is again independent of $R$ and $x_{0}$ under our assumptions. We get
$$
\|v \|_{L^\infty(B_{R}(x_{0}))} \le C_{S} R^{-N/p} \| v\|_{L^p(B_{2R}(x_{0}))}
$$ 
Serrin's Theorem (cf. Theorem 1 on page 256 of \cite{serrin}) also gives the estimate
$$
\left\| \nabla u \right\|_{L^p(B_{R}(x_{0}))} \le C_{S}R^{-1} \left\| u \right\|_{L^p(B_{2R}(x_{0}))}
$$
for solutions of \eqref{linear serrin}. Collecting these inequalities, we obtain
$$
\left\| \nabla u \right\|_{L^\infty(B_{R}(x_{0}))} \le C_{S}R^{-N/p-1} \left\| u \right\|_{L^p(B_{2R}(x_{0}))}.
$$
Using that $u(x)=o(\left| x \right|^{-\frac2{p_{0}-1}} )$, we obtain the desired estimate.\hfill\qed
\section{Proof of Theorem \ref{thsubcritical} : the subcritical case}
By Remark \ref{estimate of h}, since $p_{0}$ is subcritical, we have
\begin{equation}\label{integral estimate}  
\int_{\R^N}f(u)u\;dx<+\infty\quad\text{and}\quad \int_{\R^N}F(u)\;dx<+\infty.
\end{equation} 
Multiply equation \eqref{main} by $u\zeta$, where $\zeta$ is a standard cut-off i.e. $\zeta\equiv 1$ in $B_R$, $\zeta\equiv 0$ in $B_{2R}$ and $|\nabla\zeta|\le C/R$, $\left| \Delta\zeta\right|\le C/R^2$. Then integrate :
\begin{align*}
\int_{\R^N}\left| \nabla u\right|^2 \zeta\;dx + \int_{\R^N}u\nabla u\nabla\zeta\;dx =& \int_{\R^N}f(u)u\zeta\;dx.\\
\int_{\R^N}\left| \nabla u\right|^2 \zeta\;dx - \frac12\int_{\R^N}u^2\Delta\zeta\;dx=& 
\end{align*}
By Remark \ref{estimate of h}, the second term in the left-hand side of the above equality  converges to $0$ as $R\to+\infty$. Hence,  by monotone  convergence we have
\begin{equation}\label{p1}  
\int_{\R^N}\left| \nabla u\right|^2 \;dx = \int_{\R^N} uf(u)\;dx<+\infty
\end{equation} 
As in the classical Pohozaev identity, we may now multiply the equation by $x\cdot\nabla u\;\zeta$ and obtain
\begin{equation}\label{p2}  
\int_{\R^N}\left| \nabla u\right|^2 \;dx = \frac{2N}{N-2}\int_{\R^N} F(u)\;dx.
\end{equation} 
We now collect \eqref{p1} and \eqref{p2}. By assumption \eqref{global assumption subcritical}, if $u$ is not identically zero, then
\begin{multline*} 
\int_{\R^N}\left| \nabla u\right|^2 \;dx = \int_{\R^N} uf(u)\;dx\le (p_{0}+1)\int_{\R^N}F(u)\;dx < \frac{2N}{N-2}\int_{\R^N}F(u)\;dx \\= \int_{\R^N}\left| \nabla u\right|^2 \;dx,
\end{multline*} 
a contradiction.\hfill\qed
  
\section{Proof of Theorem \ref{th4} : the supercritical case}
In what follows, we prove Theorem \ref{th4} in the supercritical case i.e. when $p_{0}$ is in the range \eqref{supercritical range} and $f$ satisfies \eqref{lower bound}, \eqref{upper bound} and \eqref{global assumption}. 
In polar coordinates, a function $u$ takes the form $u=u(r,\sigma)$, where $r\in\R^*_{+}$, $\sigma\in S^{N-1}$, $N\ge2$, while its Laplacian is given by 
$$
\Delta u = u_{rr} + \frac{N-1}{r}u_{r} + \frac1{r^2}\Delta_{S^{N-1}}u.
$$
Recall the classical Emden change of variables and unknowns $t=\ln r$ and $u(r,\sigma) = r^{-\alpha}v(t,\sigma)$, where $\alpha=\frac2{p_{0}-1}$. Then,
\begin{align}\label{change of unknown}
v(t,\sigma)&= e^{\alpha t} u (e^t,\sigma),\nonumber\\
v_{t}(t,\sigma) &= e^{\alpha t}\left(e^t u_{r} + \alpha u\right)=e^{(\alpha+1)t}u_{r} +e^{\alpha t}\alpha u,\\
v_{tt}(t,\sigma) &= e^{\alpha t}\left(e^{2t} u_{rr} + (2\alpha+1) e^t u_{r} +\alpha^2u\right)\nonumber\\
\Delta_{S^{N-1}}v &= e^{\alpha t}\Delta_{S^{N-1}}u \nonumber
\end{align}
Writing 
\begin{equation}\label{constants ab}  
\alpha=\frac2{p_{0}-1},\quad A= \left(N-2-2\alpha\right),\quad B=\alpha^2+\alpha A,
\end{equation} 
we obtain
\begin{align*}
v_{tt} + A v_{t}&= e^{(\alpha+2)t}\left(u_{rr} + \frac{N-1}r u_{r}\right) + B e^{\alpha t}u\\
&=  e^{(\alpha+2)t}\left(-e^{-2t}\Delta_{S^{N-1}}u -f(u)\right) + B e^{\alpha t}u\\
&=  -e^{(\alpha+2)t}\left(f(e^{-\alpha t}v)\right)  - e^{\alpha t}\Delta_{S^{N-1}}v+ Bv.
\end{align*}
To summarize, $v$ solves
\begin{equation} 
v_{tt} + A v_{t} + Bv + \Delta_{S^{N-1}}v + f(e^{-\alpha t}v)e^{(\alpha+2)t} = 0\quad\text{ for $t\in\R$, $\sigma\in S^{N-1}$.}
\end{equation} 
Multiply the above equation by $v_{t}$ and integrate over $S^{N-1}$. For $t\in\R$, we find
\begin{multline} \label{energy estimate} 
\int_{S^{N-1}} \left(\frac{v_{t}^2}2\right)_{t}\;d\sigma + A \int_{S^{N-1}} v_{t}^2\;d\sigma
+ B \int_{S^{N-1}}\left(\frac{v^2}2\right)_{t}\;d\sigma\\ - \int_{S^{N-1}} \left(\frac{\left| \nabla_{S^{N-1}}v \right|^2}2\right)_{t}\;d\sigma + \int_{S^{N-1}} f(ve^{-\alpha t})v_{t}e^{(\alpha+2)t}\;d\sigma =0
\end{multline}
Let $F$ denote the antiderivative of $f$ such that $F(0)=0$. Then,
\begin{multline*}
\frac d{dt}\left[F(ve^{-\alpha t})e^{(p_{0}+1)\alpha t}\right] =\\
f(ve^{-\alpha t})\left(v_t e^{-\alpha t} -\alpha ve^{-\alpha t}\right)e^{(p_{0}+1)\alpha t} + F(ve^{-\alpha t})\alpha(p_{0}+1)e^{(p_{0}+1)\alpha t}.
\end{multline*}
So,
\begin{multline*}
f(ve^{-\alpha t})v_t e^{p_{0}\alpha t} = \\
\frac d{dt}\left[F(ve^{-\alpha t})e^{(p_{0}+1)\alpha t}\right]  + \alpha f(ve^{-\alpha t})ve^{\alpha p_{0}t} - \alpha F(ve^{-\alpha t})(p_{0}+1)e^{(p_{0}+1)\alpha t}.
\end{multline*}
Applying \eqref{global assumption},  we conclude that
$$
f(ve^{-\alpha t})v_t e^{p_{0}\alpha t} \ge \frac d{dt}\left[F(ve^{-\alpha t})e^{(p_{0}+1)\alpha t}\right] .
$$
Using this inequality in \eqref{energy estimate}, we obtain
\begin{multline*} 
\int_{S^{N-1}} \left(\frac{v_{t}^2}2\right)_{t}\;d\sigma + A \int_{S^{N-1}} v_{t}^2\;d\sigma
+ B \int_{S^{N-1}}\left(\frac{v^2}2\right)_{t}\;d\sigma\\ - \int_{S^{N-1}} \left(\frac{\left| \nabla_{S^{N-1}}v \right|^2}2\right)_{t}\;d\sigma + \int_{S^{N-1}}\frac d{dt}\left[F(ve^{-\alpha t})e^{(p_{0}+1)\alpha t}\right]\;d\sigma  \le 0.
\end{multline*}
 Integrating for $t\in (-s,s)$, $s>0$, we then derive 
\begin{multline} 
\frac12\left[\int_{S^{N-1}} v_{t}^2\;d\sigma \right]_{t=-s}^{t=s}+ A \int_{t=-s}^{t=s}\int_{S^{N-1}} v_{t}^2\;d\sigma\;dt
+ \frac B2 \left[\int_{S^{N-1}}v^2\;d\sigma\right]_{t=-s}^{t=s}\\ 
- \frac12\left[\int_{S^{N-1}} \left| \nabla_{S^{N-1}}v \right|^2\;d\sigma\right]_{t=-s}^{t=s} + \left[\int_{S^{N-1}}F(ve^{-\alpha t})e^{(p_{0}+1)\alpha t}\; d\sigma\right] _{t=-s}^{t=s}  \le 0.\label{after int}
\end{multline}
Recall the definition of $v$ given in \eqref{change of unknown} and use the improved decay estimates \eqref{fine asymptotics}: we see that $v(t,\cdot), v_{t}(t,\cdot), \left| \nabla_{S^{N-1}}v(t,\cdot)\right| $ converge to $0$ as $t\to\pm\infty$, uniformly in $\sigma\in S^{N-1}$. Passing to the limit as $s\to+\infty$ in \eqref{after int}, we finally obtain   
\begin{multline}\label{final a} 
 A \int_{\R}\int_{S^{N-1}} v_{t}^2\;d\sigma\;dt 
+\limsup_{s\to+\infty}\int_{S^{N-1}}F(ve^{-\alpha s})e^{(p_{0}+1)\alpha s}\; d\sigma \le 0.
\end{multline}
Since $p_{0}>\frac{N+2}{N-2}$, it follows from \eqref{constants ab} that $A>0$. So, both terms in \eqref{final a} are nonnegative. In particular, $v_{t}\equiv 0$ and $v$ is a function depending only on $\sigma$. Since $\lim_{t\to+\infty}v(t,\sigma)=0$ by \eqref{fine asymptotics}, we deduce that $v\equiv0$ and $u\equiv0$ as claimed.

\hfill\qed 
\noindent\bf Acknowledgments. \rm The authors wish to thank J. D\'avila for stimulating discussions on the subject. L.D. acknowledges partial support from Fondation Math\'ematique de Paris and Institut Henri Poincar\'e, where part of this work was completed.

\end{document}